\documentclass[11pt]{elsart3-1} 
\usepackage{amssymb}
\usepackage[english,francais]{babel}
\usepackage{mathabx}
 
\newtheorem{theorem}{Theorem}[section]
\newtheorem{lemma}[theorem]{Lemma}
\newtheorem{e-proposition}[theorem]{Proposition} 

\newtheorem{e-definition}[theorem]{Definition\rm}

\renewcommand{\theequation}{\arabic{equation}}
\def\og{\leavevmode\raise.3ex\hbox{$\scriptscriptstyle\langle\!\langle$~}}
\def\fg{\leavevmode\raise.3ex\hbox{~$\!\scriptscriptstyle\,\rangle\!\rangle$}}

\newcommand \newperp  {\underline{\partial}_\perp}
\newcommand \sbar {{\bar s}}
\newcommand \bel {\be \label}

\newcommand \hb{\overline h}

\newcommand \Kcal {\mathcal K}
\newcommand \Hcal {\mathcal H}
\newcommand \Boxt {\widetilde {\Box}}
\newcommand \del \partial
\newcommand \delu {\underline{\del}}

\newcommand \RR{\mathbb{R}}
\newcommand {\vep}{\varepsilon}

\let\oldmarginpar\marginpar
\renewcommand\marginpar[1]{\-\oldmarginpar[\raggedleft\footnotesize #1]%
{\raggedright\footnotesize #1}}

\newcommand \Hf H 
\newcommand \be {\begin{equation}}
\newcommand \ee {\end{equation}}

\newcommand \bei {\begin{itemize}}
\newcommand \eei {\end{itemize}}

\journal{the Acad\'emie des sciences}
\begin{document}
% place in the next line the header (rubrique) chosen for your article,
% if you know it (you can also have 2, format : Header1/Header2
\centerline{}
\begin{frontmatter}

% Title, authors and addresses

% use the thanksref command within \title, \author or \address for footnotes;
% use the ead command for the email address,
% and the form \ead[url] for the home page:
% \title{Title\thanksref{label1}}
% \thanks[label1]{}
% \author{Name\thanksref{label2}}
% \ead{email address}
% \ead[url]{home page}
% \thanks[label2]{}
% \address{Address\thanksref{label3}}
% \thanks[label3]{}
\selectlanguage{english}
\title{The global nonlinear stability of Minkowski space
\\
for the Einstein equations in presence of a massive field}

% use optional labels to link authors explicitly to addresses:
% \author[label1,label2]{}
% \address[label1]{}
% \address[label2]{}
% The [label1] can be suppressed if there is only one address for all authors

\selectlanguage{english}
\author[authorlabel1]{Philippe G. LeFloch}
\ead{ contact@philippelefloch.org}
\and 
\author[authorlabel2]{Yue Ma}
\ead{yuemath@mail.xjtu.edu.cn}

\address[authorlabel1]{Laboratoire Jacques-Louis Lions, Centre National de la Recherche Scientifique,
\\
Universit\'e Pierre et Marie Curie, 4 Place Jussieu, 75252 Paris, France. }

\address[authorlabel2]{School of Mathematics and Statistics, Xi'an Jiaotong University,
\\
 Xi'an, 710049 Shaanxi, Popular Republic of China.}

% If you know the dates of reception, and acceptation you can put them now;
%  idem the name of the person presenting the Note

\medskip
%\begin{center}
%{\small Received *****; accepted after revision +++++\\
%Presented by}
%\end{center}

\begin{abstract}
\selectlanguage{english} We provide a significant extension of the Hyperboloidal Foliation Method introduced by the authors in 2014 in order to establish global existence results for systems of quasilinear wave equations posed on a curved space, when wave equations and Klein-Gordon equations are coupled. 
This method is based on a $(3+1)$ foliation (of the interior of a future light cone in Minkowski spacetime) by spacelike and asymptotically hyperboloidal hypersurfaces. In the new formulation of the method, we 
succeed to cover wave-Klein-Gordon systems containing ``strong interaction'' terms at the level of the metric, and then generalize our method in order to esstablish a new existence theory for the Einstein equations of general relativity. 
Following pioneering work by Lindblad and Rodnianski on the Einstein equations in wave coordinates, we establish the nonlinear stability of Minkowski spacetime for self-gravitating massive scalar fields.

% {\it To cite this article: P.G. LeFloch \& Y. Ma, C. R. Acad. Sci. Paris, Ser. I (2016).}

\vskip 0.5\baselineskip

\

{\sl \large Stabilit\'e nonlin\'eaire globale de l'espace-temps de Minkowski pour les champs massifs. }

\

\selectlanguage{francais}
% Text of abstract in French
\noindent{\bf R\'esum\'e} \vskip 0.5\baselineskip \noindent

Nous g\'en\'eralisons la M\'ethode du Feuiletage Hyperboloidal introduite par les auteurs en 2014
pour traiter des syst\`emes quasilin\'eaires couplant des \'equations d'ondes et des \'equations de Klein-Gordon. Dans cette nouvelle formulation, nous r\'eussissons \`a traiter des termes m\'etriques ``d'interaction forte''.  Nous appliquons cette m\'ethode pour d\'emontrer la stabilit\'e nonlin\'eaire de l'espace de Minkowski pour les \'equations d'Einstein des champs scalaires massifs auto-gravitants. En suivant un travail de Lindblad et de Rodnianski, nous analysons la structure des \'equations d'Einstein en coordonn\'ees d'ondes, qui constituent pr\'ecis\'ement un syst\`eme d'\'equations d'ondes quasi-lin\'eaires avec ``interaction forte''. 

%{\it Pour citer cet article~: P.G. LeFloch \& Y. Ma, C. R. Acad. Sci. Paris, Ser. I (2016).}

\end{abstract}
\end{frontmatter}

% now the Version fran�aise abr�g�e, if it exists 
\selectlanguage{francais}
\section*{Version fran\c{c}aise abr\'eg\'ee}
% Text of your Version fran�aise abr�g�e here.
% Note you do not need to repeat here equations that you use in the
% main text - for example 'voir (3)' is quite acceptable.

Nous \'etudions le probl\`eme de l'existence globale en temps de solutions r\'eguli\`eres d'\'equations d'ondes nonlin\'eaires, avec deux objectifs principaux: 

\bei 

\item Nous g\'en\'eralisons la m\'ethode du feuiletage hyperboloidal \cite{PLF-YM-book} propos\'ee par les auteurs 
en 2014.  

\item Cette m\'ethode nous permet d'analyser les \'equations d'Einstein et de d\'emontrer la stabilit\'e nonlin\'eaire de l'espace de Minkowski en pr\'esence d'un champ scalaire massif auto-gravitant. 

\eei 

Les aspects nouveaux de notre m\'ethode (voir le texte en anglais pour plus de d\'etails) sont les suivants:

\bei 

\item[(1)]  inegalit\'es de Sobolev et de Hardy associ\'ees au feuilletage hyperboloidal de l'int\'erieur d'un c\^one de lumi\`ere (de l'espace-temps de Minkowski). 

\item[(2)] Estimations de type $L^\infty$--$L^\infty$ pour l'\'equation d'onde et pour l'\'equation de Klein-Gordon sur un espace courbe.  

\item[(3)] Une hi\'erarchie d'estimations d'\'energie ayant des  croissances alg\'ebriques en temps diff\'erentes.

\eei 

\

Par ailleurs, en ce qui concerne le traitement des \'equations d'Einstein proprement d\^{\i}t, nous combinons diff\'erentes id\'ees et techniques, \'enonc\'ees comme suit: 

\bei 
 
\item[(4)] d\'ecomposition de la courbure de Ricci (``null forms'', ``quasi-null forms''), 

\item[(5)] d\'ecomposition de tenseurs bas\'ee sur la condition d'onde, 

\item[(6)] int\'egration le long de caracteristiques, 
et 

\item[(7)] une hierarchie d'\'energies adapt\'ees \`a la structure des  \'equations d'Einstein-Klein-Gordon.  

\eei 

%===================================================================

\selectlanguage{english} 

\section*{English version}
\section{Introduction}
\label{}

We are interested in the global-in-time existence problem for small amplitude solutions to nonlinear wave equations posed on a curved spacetime, with a two-fold objective. First, we provide a significant extension of the Hyperboloidal Foliation Method, proposed by the authors \cite{PLF-YM-book} in 2014. This method is based on a $3+1$ foliation (of the interior of a future light cone in Minkowski spacetime) by spacelike hyperboloidal hypersurfaces and on Sobolev and Hardy-type inequalities adapted to this hyperboloidal foliation.  This method applies to a broad class of coupled nonlinear wave-Klein-Gordon systems on curved space, and takes its root in a work by Klainerman \cite{Klainerman85} and H\"ormander \cite{Hormander} on the Klein-Gordon equation.  

 In comparison to our earlier formulation  \cite{PLF-YM-book}, we are now able to encompass a broader class of coupled systems involving ``strong interaction'' terms (as we call them, see below). Recall that Klainerman used the decomposition of the wave operator along an hyperboloid al foliation and was able to establish a (non-sharp) decay rate of $t^{-5/4}$ and conclude with a theory of existence of quasilinear Klein-Gordon equations.  
Later on, Hormander worked directly with the energy on hyperboloidal hypersurfaces and Sobolev inequalities and derived the sharp rate $t^{-3/2}$ for the problem posed in flat space. 

Our work \cite{PLF-YM-book} provides a way to combine both approaches and encompass {\sl systems of coupled} wave and Klein-Gordon equations. Importantly, we work directly {\sl within the hyperboloidal foliation}  and, in order to encompass equations posed on curved space, we establish sup-norm bounds leading us to the sharp rate $t^{-3/2}$.  

Our second objective is to apply this method to the {\sl Einstein equations of general relativity} and, by constructing spacelike and asymptotically hyperboloidal hypersurfaces, 
we offer a new strategy of proof in order to establish the nonlinear stability of Minkowski spacetime. Our method applies {\sl self-gravitating massive} scalar fields, while all earlier works were restricted to vacuum spacetimes or spacetimes with massless scalar fields; cf.~the pioneering work by Christodoulou and Klainerman \cite{CK}, the proof by Lindblad and Rodnianski \cite{LR1,LR2} (based on the wave gauge), and the extension by Bieri and Zipser \cite{BieriZipser}.

One of the simplest  wave-Klein-Gordon model is, in flat space, $\Box u = P(\del u, \del v)$, 
$\Box v + v = Q(\del u, \del v)$, where $P, Q$ are quadratic forms in the first-order derivatives $\del u = (\del_\alpha u)$ and $\del v = (\del_\alpha v)$ and the two unknowns $u, v$ are defined over Minkowski space $\RR^{3+1}$. (Here $\alpha=0,1, 2, 3$.)
Many models arising in mathematical physics involve interactions between massive and massless fields. Let us mention the Dirac-Klein-Gordon equations, the Proca equation (massive spin-$1$ field in Minkowski spacetime), the Einstein-massive field system, and the field equations of modified gravity described by the Hilbert-Einstein functional $\int_M f(R_g) \, dv_g$ with, typically,  $f(R_g) = R_g + \kappa (R_g)^2$ and $\kappa>0$, where $R_g$ is the scalar curvature of a Lorentzian manifold $(M,g)$.

The vector field method was introduced by Klainerman \cite{Klainerman80,Klainerman85} around 1980--1985. It primarily applies to quasilinear wave equations posed on the $(3+1)$-dimensional  Minkowski spacetime, and leads to 
global-in-time well-posedness results when the initial data are sufficiently small in some Sobolev spaces. 
The method relies on the use of the conformal Killing fields of Minkowski space, 
suitably weighted energy estimates, and the so-called Klainerman-Sobolev inequalities. 
Nonlinearities are assumed to satisfy the `null condition' and a bootstrap argument is formulated  and relies on time decay estimates. In comparison, quasilinear Klein-Gordon equations have atracted less attention in the literature, despite pioneering contributions by Klainerman \cite{Klainerman85}, Shatah \cite{Shatah}, and H\"ormander \cite{Hormander}.  
 
In our work \cite{PLF-YM-book,PLF-YM-zero,PLF-YM-un,PLF-YM-deux}, we have addressed this major challenge of developing a method for quasilinear wave-Klein-Gordon systems on curved space. The main difficulty comes from the fact that a smaller symmetry group is available to deal with Klein-Gordon equations, since the scaling field $t\del_t + r\del_r$ does not commute with the linear Klein-Gordon operator in flat spacetime. 
While additional decay for Klein-Gordon equations, that is, $t^{-3/2}$ in four dimensions is available 
(solutions to wave equations decaying only like $t^{-1}$), a robust technique to deal with the {\sl coupling} of wave equations and Klein-Gordon equations is required and we have developped the Hyperboloidal Foliation Method precisely for that purpose. 
For earlier contributions on the analysis of Klein-Gordon equation with a limited number of Killing fields, we refer to Katayama \cite{Katayama12a} and the references cited therein. 

 %==============================================================================
 
\section{The hyperboloidal foliation method} 
  
Following \cite{PLF-YM-book}, we rely solely on the Lorentz boosts (or hyperbolic rotations), which generate a foliation (of the interior of a light cone) of Minkowski spacetime by hyperboloidal hypersurfaces (that is, surfaces of constant distance from a base point). We introduce here Lorentz-invariant energy norms based on these boosts and are able to revisite all the standard arguments of the so-called vector field method and carefully analyze the energy flux on the hyperboloids. With a suitable extensions of Sobolev and Hardy inequalities adapted to the hyperboloidal foliation, we are able to encompass a broad class of coupled systems. 

The method in \cite{PLF-YM-book} is based on the {\sl hyperboloidal hypersurfaces}
$\Hcal_s := \big\{ (t,x) \, \big/ \, t>0; \,  t^2 - |x|^2 =s^2 \big\}$
 parametrized by their hyperbolic radius $s>s_0>1$, 
and consider the  {\sl foliation of the future light cone} $\Kcal := \big\{ (t,x) \, / \, |x|\leq t-1 \big\}$. Note in passing that $s \leq t \leq s^2$.  
We imposed initial data on the hypersurface $t=s_0>1$ or directly on the hyperboloid $s=s_0$. Our 
energy estimates therein are formulated in domains limited by two hyperboloids, that is,  
$\begin{array}{ccc}
\Kcal_{[s_0,s_1]} 
:=& \big\{(t,x) \, / \,  |x| < t-1,\, (s_0)^2 \leq t^2- |x|^2\leq (s_1)^2,{t>0} \big\}. 
\end{array}
$
 Our analysis is performed in the {\sl semi-hyperboloidal frame} (as we propose to call it), consisting of the Lorentz boosts
$L_a := x_a \del_t  + t \del_a, \, a=1,2,3$ 
and a time-like vector. More precisely, by definition, 
this frame consists of the following three vectors tangent to the hyperboloids 
$
\delu_a := {L_a \over t}
$
and the timelike vector 
$
\delu_0 := \del_t. 
$
Accordingly, we have the semi-hyperboloidal decomposition of the (flat) wave operator: 
$
\Box u = -\frac{s^2}{t^2} \delu_0 \delu_0 u - \frac{x^a}{t} \delu_0 \delu_a u - \frac{x^a}{t}\delu_a \delu_0 u + \sum_a\delu_a\delu_a u - \frac{3}{t}\del_tu$.   
(In comparison, the standard choice in the literature is the `null frame', containing three vectors tangent to the light cone.)

The {\sl hyperboloidal energy} associated with the hypersurface $\Hcal_s$ involves certain weighted derivatives, 
and we want to point out that we will use the full expression of the corresponding energy flux on the hyperboloids. 
Let us also mention one important functional inequality.  

\

\begin{lemma}[Sobolev estimate on hyperboloids]
For all functions $u$ defined on a hyperboloid $\Hcal_s$ in Min\-kowski space $\RR^{3+1}$ and with sufficiently fast decay, one has 
$
\sup_{(t,x) \in \Hcal_s} t^{3/2} |u(t,x)| \lesssim \sum_{|I|\leq 2} \| L^I u \|_{L^2(\Hcal_s)}$
(for $s \geq s_0 >1$) 
with summation over $L \in \big\{ L_a = x_a \del_t  + t \del_a \big\}$, where $I$ denotes a multi-index. 
\end{lemma}
  
\

The hyperboloidal foliation method is based on a hierarchy of bounds for the curved metric and the source-terms.
Specifically, in \cite{PLF-YM-book}, we used three levels of regularity and algebraic growth rates and, remarkably, our bound is uniform for the low-order energy of wave components. 

%====================================================================== 

\section{The Einstein-massive field system} 

We present a new method for proving the nonlinear stability of Minkowski spacetime, which applies to self-gravitating massive scalar fields \cite{PLF-YM-un,PLF-YM-deux}. The statement of the problem is as follows (following Choquet-Bruhat et al.  \cite{CB}): 
we search for a spacetime $(M,g)$ satisfying the {\sl Einstein equations} 
$$
R_{\alpha\beta} - {R \over 2} g_{\alpha\beta} = T_{\alpha\beta}
$$
for the stress-energy tensor of a scalar field $\phi$, that is,  
$$
T_{\alpha\beta} := \nabla_\alpha \phi \nabla_\beta \phi - \Big( {1 \over 2} \nabla_\gamma \phi \nabla^\gamma \phi + V(\phi) \Big) g_{\alpha\beta}, 
$$
where the potential is taken to be $V(\phi) := \frac{c^2}{2}\phi^2$ ($c>0$ being the mass of the scalar field). Using the contracted Bianchi identities, it is not difficult to derive the {\sl Klein-Gordon equation} $\Box_g\phi = V'(\phi) = c^2 \phi$. 
Our objective is to study the associated Cauchy problem when the initial data set is a perturbation of a spacelike hypersurface in Minkowski space. 

\

\begin{theorem}[Nonlinear stability of Minkowski space for massive fields]
Consider the Einstein-scalar field system in wave coordinates (that is, $\Box_g x^\alpha  = 0$): 
$$
\begin{array}{ccc}
 \Boxt_g g_{\alpha\beta}  = \big( Q_{\alpha\beta} + P_{\alpha\beta}\big)(g;\del g,\del g)
- 2 \big(\del_\alpha\phi\del_\beta\phi + V(\phi)g_{\alpha\beta}\big),  
\qquad
\Boxt_g \phi - V'(\phi) = 0,  
\end{array}
$$
where $\Boxt_g := g^{\alpha'\beta'}\del_{\alpha'}\del_{\beta'}$ is the so-called reduced wave operator,
$Q_{\alpha\beta}$ are null terms, and $P_{\alpha\beta}$ are ``weak null'' terms.
Consider an initial data set $(\overline M, \overline g, k)$ which is close to a spacelike slice of Minkowski space and is asymptotically hyperboloidal (in a suitable sense) 
Then, the initial value problem for the Einstein-massive field system admits a global solution in wave coordinates, which defines a future geodesically complete spacetime $(M,g)$.
\end{theorem} 

\

In \cite{PLF-YM-un,PLF-YM-deux}, the massive field is assumed to be compactly supported, but this assumption is not essential and is removed in \cite{PLF-YM-trois}, which also treats the more general theory of modified gravity.  
Our proof relies on the wave gauge, after the pioneering work of Lindblad and Rodnianski \cite{LR1,LR2}.  The following challenges and techniques are in order: 

\bei 

\item {\sl Tensorial structure}: The geometric structure of the Einstein equations in combination with the wave coordinate condition $\Box_g x^\alpha = 0$ allows us to decompose the quadratic nonlinearities as a sum of null terms and ``weak null'' terms. 

\item {\sl Hyperboloidal foliation}: Having fewer Killing fields at our disposal, we rely on the foliation generated by the Lorentz boosts, that is, the hyperboloids of Minkowski space and we introduce Lorentz-invariant energy norms. As explained in Section~2, we need to establish Sobolev and Hardy-type inequalities on hyperboloids. 

\item {\sl Sharp pointwise estimates:} We derive $L^\infty$--$L^\infty$ estimates for, both, the wave equation and  Klein-Gordon equations on curved space. We use a technique of integration along well-chosen curves (see below). 

\item {\sl Hierarchy of energy bounds}: Several levels of regularity and time growth rates 
are required in  our bootstrap argument, and successive improvements of the estimates are performed in the proof. 

\eei 
In the rest of this text, we present some of our technique for a simplified, but challenging, model. 

\section{Wave-Klein-Gordon model with strong interactions}

Consider the following system 
with strong interactions at the metric level, 
which (we formally extract from the Einstein equations in wave coordinates): 
\bel{eq:MODP}
-\Box u = P^{\alpha\beta}\del_\alpha v\del_\beta v + Rv^2, 
\qquad -\Box v + u \, H^{\alpha\beta} \del_\alpha\del_\beta v + c^2 v = 0. 
\ee

\begin{theorem}[Nonlinear wave-Klein-Gordon model with strong interaction]
Consider the nonlinear wave-Klein-Gordon model (\ref{eq:MODP}) 
with given constants $P^{\alpha\beta}, R, H^{\alpha\beta}$ and $c>0$. 
For any $N \geq 8$, there exists $\epsilon=\epsilon(N)>0$ such that  if the initial data satisfy 
$
\| (u_0, v_0) \|_{H^{N+1}(\RR^3)} + \| (u_1, v_1) \|_{H^N(\RR^3)} < \epsilon
$
 then the Cauchy problem for (\ref{eq:MODP}) admits a global-in-time solution.
\end{theorem}
 
\

We proceed with a bootstrap argument based on a hierarchy of energy bounds posed along the hyperboloidal foliation (stated in (\ref{ineq energy assumption}), below). The proof of Proposition~\ref{propo44}  relies in particular on a suitable  decomposition of the Klein Gordon equation on curved space as well as an involved ODE lemma. 
 
  \

\begin{e-proposition}[$L^\infty$--$L^\infty$ estimate for Klein-Gordon equations on curved space]
\label{propo44} 
Consider the Klein-Gordon equation on a curved background $-\Boxt_g v + c^2 v = f$ with metric $g^{\alpha\beta} = m^{\alpha\beta} - h^{\alpha\beta}$ given by a perturbation of the Minkowski metric, and with compactly supported data prescribed on a hyperboloid 
$v|_{\Hcal_{s_0}} = v_0, \,  \del_t v|_{\Hcal_{s_0}} = v_1$
 for sufficiently smooth and spatially compactly supported data $v_0, v_1$. 
Then, in the future of $\Hcal_{s_0}$, one has 
$$
s^{3/2}|v(t,x)| + {t \over s} \, s^{3/2}|\newperp  v(t,x)|\lesssim V(t,x)
$$
with $V$ defined below and $\newperp := \del_t + \frac{x^a}{t}\del_a$. 
\end{e-proposition}

\

Note that $\newperp$ is orthogonal to the hyperboloids for the Minkowski metric and  coincides, up to an essential factor $1/t$, with the scaling vector field. 
We use the notation 
$h_{t,x}(\lambda):= \hb^{00}(\lambda t/s,\lambda x/s)$ (with $ s^2 = t^2 - r^2$) and 
consider the derivative in $\lambda$, that is 
$$
h_{t,x}'(\lambda) 
 = {t \over s} \del_t \hb^{00}(\lambda t/s,\lambda x/s) + {x^a \over s} \del_a \hb^{00}(\lambda t/s,\lambda x/s)
= {t \over s} \, \newperp \hb^{00}(\lambda t/s,\lambda x/s).
$$
Fix a constant $C>0$ (chosen later on) and define the function $V$ first  {\sl ``far" from the light cone} $0\leq r/t\leq {s_0^2 - 1 \over 1+s_0^2}$: 
$$
\begin{array}{ccc}
V(t,x) :=  \, & \Big( \|v_0\|_{L^\infty(\Hcal_{s_0})} + \|v_1\|_{L^\infty(\Hcal_{s_0})} \Big)
\Big(1+\int_{s_0}^s | h_{t,x}'(\sbar)|e^{C\int_\sbar^s|h_{t,x}'(\lambda)|d\lambda} \, d\sbar \Big)
\\
& + F(s) + \int_{s_0}^s F(\sbar)|h_{t,x}'(\lambda)|e^{C\int_\sbar^s|h_{t,x}'(\lambda)|d\lambda} \, d\sbar
\end{array}
$$
and then {\sl ``near'' the light cone}  ${s_0^2 - 1 \over 1+s_0^2} <r/t<1$ by 
$$
V(t,x) :=
F(s) + \int_{S(r/t)}^s F(\sbar)|h_{t,x}'(\sbar)|e^{C\int_\sbar^s|h_{t,x}'(\lambda)|d\lambda} \, d\sbar 
$$
with $S(r/t) :=  \sqrt{\frac{t+r}{t-r}}$. Here, $F$ is defined by a suitable integration of the given source-term $f$.
 
 \

\begin{e-proposition}[$L^\infty$--$L^\infty$ estimate for the wave equation with source]
\label{porpo45}
Let $u$ be a spatially compactly supported solution to the wave equation $-\Box u = f$ with vanishing initial data 
and source $f$ satisfying  
$
|f|\lesssim {1 \over t^{2+\nu} (t-r)^{1-\mu}}$ ($t \geq 2$) 
for some exponents $0<\mu\leq 1/2$ and $0< |\nu|\leq 1/2$: 
$$ 
|u(t,x)|
\lesssim \left\{ 
\begin{array} {cc} 
{1 \over \nu\mu} {1 \over (t-r)^{\nu-\mu} \, t}, \qquad & 0< \nu\leq 1/2, 
\\ 
{1 \over |\nu|\mu} {(t-r)^{\mu} \over t^{1+\nu}}, &-1/2\leq \nu < 0.  
\end{array}\right.
$$
\end{e-proposition}

\

The proof is based on the explicit representation formula for the wave equation. We now sketch our bootstrap argument, based on the following assumptions (with $|J|=k$): : 
\bel{ineq energy assumption}
\begin{array} {ccccc}
&E(s,\del^IL^J u)^{1/2} + s^{-1/2} E(s,\del^IL^J v)^{1/2} \leq C_1\vep s^{k\delta}, & \quad &|I|+|J|\leq N, 
\\
&E(s,\del^IL^J u)^{1/2}\leq C_1\vep, &\quad &|I|+|J|\leq N-4,   
\\
&E(s,\del^IL^J v)^{1/2}\leq C_1\vep s^{k\delta}, &\quad &|I|+|J|\leq N-4. 
\end{array} 
\ee
Using Sobolev and Hardy inequalities adapted to the hyperboloids,   
we then deduce basic decay bounds. 
Next, we derive the following sharp sup-norm bounds for $|I|+|J|\leq N-4$, which are at the heart of our argument: 
\be
\begin{array}{ccc}
\sup_{\Hcal_s} t|L^J u|
+
\sup_{\Hcal_s} \, (t/s)^{1/2 -4\delta}t^{3/2}|\del^IL^J v| 
+
\sup_{\Hcal_s} \, (t/s)^{3/2-4\delta}t^{3/2}|\newperp \del^IL^J v| 
&\lesssim C_1\vep s^{k\delta}. 
\end{array}
\ee
We proceed 
as follows: 
\begin{itemize}

\item {\sl -- First bound for the wave component} ($L^\infty$--$L^\infty$ estimate for wave equations) for $|I|+|J|\leq N-7$:   
$$
| \del^IL^J u | \lesssim C_1 \vep t^{-3/2} + (C_1\vep)^2 (t/s)^{-(k+4)\delta}t^{-1}s^{(k+4)\delta}.
$$

\item 
{\sl Second bound for the wave component and  
 first bound for the Klein-Gordon component} ($L^\infty$--$L^\infty$ for wave and K-G equations))
$$
\begin{array}{ccc}
|u(t,x)| \lesssim C_1\vep t^{-1}; \qquad 
|v| +
{t \over s} |\newperp  v(t,x)|
\lesssim C_1\vep (t/s)^{-2+7\delta}s^{-3/2}.
\end{array}
$$ 

\item 
{\sl Second bound for the Klein-Gordon component} (again the $L^\infty$--$L^\infty$ for K-G) 
  for $|I|\leq N-4$: 
$$
|\newperp  \del^I v(t,x)|\lesssim C_1\vep(t/s)^{-3/2+4\delta} t^{-3/2}; 
\qquad 
|\del^I v(t,x)|\lesssim C_1 \vep (t/s)^{-1/2+4\delta}t^{-3/2}.
$$ 

\item 
{\sl Third {\rm (and sharp, except for the higher order derivatives of $v$)} \sl bound for the wave and Klein-Gordon components} for  $|I|+|J|\leq N-4$: 
$$
\sup_{\Hcal_s}\big(t|L^J u|\big)\lesssim C_1\vep s^{k\delta},
$$
$$
\sup_{\Hcal_s}\big((t/s)^{3 - 7\delta}s^{3/2}|\newperp \del^IL^J v|\big) 
+ \sup_{\Hcal_s}\big((t/s)^{2-7\delta}s^{3/2}|\del^IL^J v| \big)\lesssim C_1\vep s^{k\delta},
$$
$$
\sup_{\Hcal_s}\big((t/s)^{1-7\delta}s^{3/2}|\del_\alpha \del^IL^J v|\big)\lesssim C_1\vep s^{k\delta}.
$$
\end{itemize} 
The pointwise estimates follow from Propositions \ref{propo44} and \ref{porpo45}, the bootstrap assumptions, the structure of the equations, and various commutation properties enjoyed by the vector fields under consideration in the Hyperboloidal Foliaiton Method. 
Finally,we can conclude and close our bootstrap argument by returning to the (differentiated) system 
$$ 
-\Box \del^I L^Ju 
 = \del^IL^J\left(P^{\alpha\beta}\del_{\alpha}v\del_{\beta}v\right) + \del^IL^J\left(Rv^2\right), 
$$
$$
-\Box \del^IL^J v + u \, H^{\alpha\beta}\del^IL^J v+ c^2\del^IL^J v 
 = -[\del^IL^J, u \, H^{\alpha\beta}\del_{\alpha}\del_{\beta}]v,
$$ 
and showing that all source-terms provide integrable contributions to the energy.

 %=============================================================== 


\begin{thebibliography}{00} 

\bibitem{BieriZipser} {\sc L. Bieri and N. Zipser,}
{\sl Extensions of the stability theorem of the Minkowski space in general relativity,}
AMS/IP Studies Adv. Math. 45. Amer. Math. Soc., International Press, Cambridge,  2009.

\bibitem{CB} {\sc Y. Choquet-Bruhat},
{\sl General relativity and the Einstein equations}, Oxford Math. Monograph,
Oxford University Press, 2009.

\bibitem{CK} {\sc D. Christodoulou and S. Klainerman,}
{\sl The global nonlinear stability of the Minkowski space,}
Princeton Math. Ser. 41, 1993.  

\bibitem{Hormander}{\sc L. H\"ormander,}
{\sl Lectures on nonlinear hyperbolic differential equations,}
Springer Verlag, Berlin, 1997.

\bibitem{Katayama12a} {\sc S. Katayama,}
Global existence for coupled systems of nonlinear wave and Klein-Gordon equations in three space dimensions,
Math. Z. 270 (2012), 487--513. 

\bibitem{Klainerman80}{\sc S. Klainerman,}
Global existence for nonlinear wave equations,
 Comm. Pure Appl. Math. 33 (1980), 43--101.

\bibitem{Klainerman85}{\sc S. Klainerman,}
Global existence of small amplitude solutions to nonlinear Klein-Gordon equations in four spacetime dimensions,
Comm. Pure Appl. Math. 38 (1985), 631--641.

\bibitem{PLF-YM-book} {\sc P.G. LeFloch and Y. Ma,}
{\sl The hyperboloidal foliation method,}  World Scientific Press, Singapore, 2014.

\bibitem{PLF-YM-zero}{\sc P.G. LeFloch and Y. Ma},
The mathematical validity of the f(R) theory of modified gravity, ArXiv:1412.8151. 

\bibitem{PLF-YM-un}{\sc P.G. LeFloch and Y. Ma},
The global nonlinear stability of Minkowski space for self-gravitating massive fields. The wave-Klein-Gordon model, Comm. Math. Phys. (2016), to appear.  

\bibitem{PLF-YM-deux}{\sc P.G. LeFloch and Y. Ma},
The global nonlinear stability of Minkowski space for self-gravitating massive fields, 
Preprint ArXiv:1511.03324.

\bibitem{PLF-YM-trois}{\sc P.G. LeFloch and Y. Ma},
The nonlinear stability of Minkowski spacetime for massive matter in Einstein's theory and f(R)-gravity, in preparation. 

\bibitem{LR1} {\sc H. Lindblad and I. Rodnianski,}
Global existence for the Einstein vacuum equations in wave coordinates,
Comm. Math. Phys. 256 (2005), 43--110.

\bibitem{LR2} {\sc H. Lindblad and I. Rodnianski,}
The global stability of Minkowski spacetime in harmonic gauge,
Ann. of Math. 171 (2010), 1401--1477.

\bibitem{Shatah} {\sc J. Shatah,}
Normal forms and quadratic nonlinear Klein-Gordon equations,
Comm. Pure Appl. Math. 38 (1985), 685--696.

\end{thebibliography}
\end{document}